\documentstyle[amscd,amssymb,verbatim,epsf]{amsart}

\begin{document}
\theoremstyle{plain}
\newtheorem{Thm}{Theorem}
\newtheorem{IMThm}{Main Theorem}
\newtheorem{Cor}{Corollary}
\newtheorem{Ex}{Example}
\newtheorem{Con}{Conjecture}
\newtheorem{Main}{Main Theorem}
\newtheorem{Lem}{Lemma}
\newtheorem{Prop}{Proposition}

\theoremstyle{definition}
\newtheorem{Def}{Definition}
\newtheorem{Note}{Note}

\theoremstyle{remark}
\newtheorem{notation}{Notation}
\renewcommand{\thenotation}{}

\errorcontextlines=0
\numberwithin{equation}{section}
\renewcommand{\rm}{\normalshape}%

\title[Neutral K\"{a}hler Metric]%
   {A Neutral K\"{a}hler Metric on the Space of Time-like Lines in Lorentzian 3-space}
\author{Brendan Guilfoyle}
\address{Brendan Guilfoyle\\
          Department of Mathematics and Computing \\
          Institute of Technology, Tralee \\
          Clash \\
          Tralee  \\
          Co. Kerry \\
          Ireland.}
\email{brendan.guilfoyle@@ittralee.ie}
\author{Wilhelm Klingenberg}
\address{Wilhelm Klingenberg\\
 Department of Mathematical Sciences\\
 University of Durham\\
 Durham DH1 3LE\\
 United Kingdom.}
\email{wilhelm.klingenberg@@durham.ac.uk }

\date{November 24th, 2006}

\begin{abstract}
We study the neutral K\"ahler metric on the space of time-like lines in Lorentzian ${\Bbb{E}}^3_1$, which we
identify with the total space of the tangent bundle to the hyperbolic plane. We find all of the infinitesimal
isometries of this metric, as well as the geodesics, and interpret them in terms of the Lorentzian metric on
${\Bbb{E}}^3_1$. In addition, we give a new characterisation of Weingarten surfaces in Euclidean 
${\Bbb{E}}^3$ and Lorentzian ${\Bbb{E}}^3_1$ as the vanishing of the scalar curvature of the associated 
normal congruence in the space of oriented lines. Finally, we relate our construction to the classical 
Weierstrass representation of minimal and maximal surfaces in ${\Bbb{E}}^3$ and ${\Bbb{E}}^3_1$.
\end{abstract}

\maketitle

\section{Time-like lines in ${\Bbb{E}}^3_1$}

Consider the Lorentzian space ${\Bbb{E}}^3_1$ with flat coordinates ($x^1,x^2,x^3$) so that the metric takes the form:
\[
ds^2=\left(dx^1\right)^2+\left(dx^2\right)^2-\left(dx^3\right)^2.
\]
For convenience we let $z=x^1+ix^2$ and $t=x^3$.

The set of future-pointing time-like geodesics in ${\Bbb{E}}^3_1$ can be identified with TH$^2$, the tangent bundle 
to the hyperbolic plane, as follows. H$^2$ is one of the two connected components of the 2-sheeted hyperboloid
in ${\Bbb{E}}^3_1$:
\[
\left(x^1\right)^2+\left(x^2\right)^2-\left(x^3\right)^2=-1.
\]
This embedded disc can be parameterized by the map $h:{\Bbb{C}}\rightarrow{\Bbb{E}}^3_1$:
\[
z=\frac{2\xi}{1-\xi\bar{\xi}} \qquad\qquad t=\frac{1+\xi\bar{\xi}}{1-\xi\bar{\xi}},
\]
for $|\xi|<1$. The metric g induced on H$^2$ by the Lorentzian metric is positive definite with 
coordinate expression:
\[
ds^2=\frac{4}{(1-\xi\bar{\xi})^2}d\xi d\bar{\xi}.
\]
In particular, the metric is of constant curvature -1 and $\xi$ is a holomorphic coordinate on H$^2$.

The space ${\Bbb{L}}^{3}_{1,-}$ of future-pointing time-like lines can, by parallel transport, be identified with the
tangent bundle to this embedded disc. The direction vector (normalised to length -1) of the line determines the point 
on the hyperboloid, while the Lorentz orthogonal vector from the line to the origin is the tangent vector to the 
disc. Thus we have identified ${\Bbb{L}}^{3}_{1,-}$ with TH$^2$ (see figure below).

\vspace{0.1in}
\setlength{\epsfxsize}{3.5in}
\begin{center}
   \mbox{\epsfbox{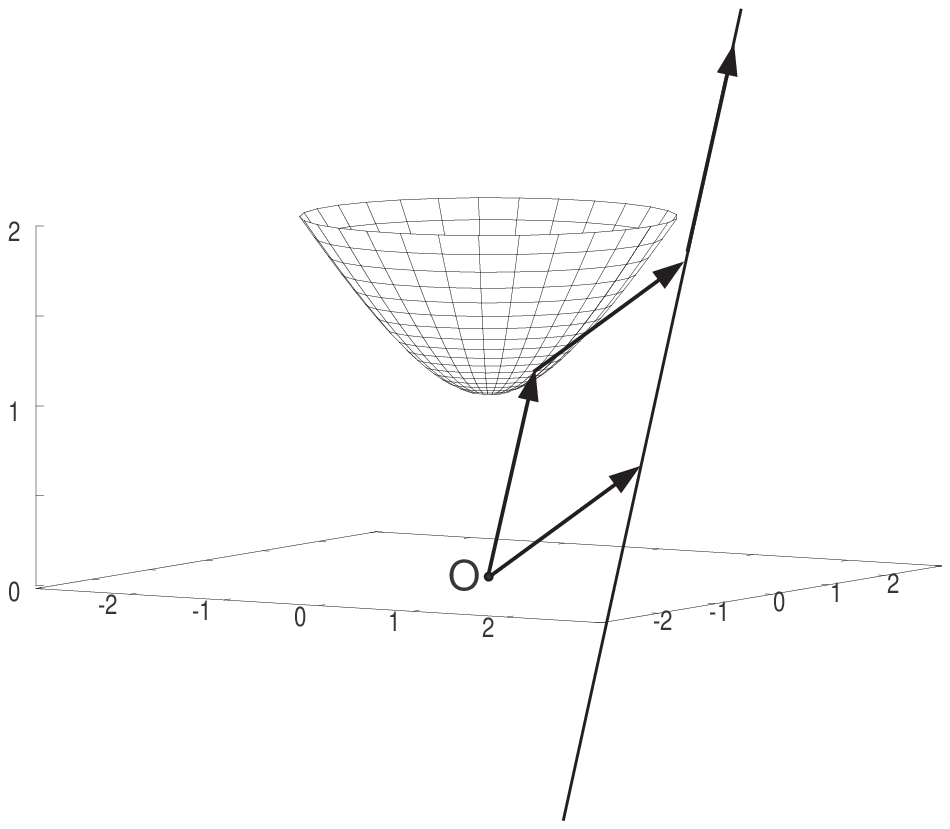}}
\end{center}
\vspace{0.1in}

In \cite{gak4} a method was presented for constructing a K\"ahler structure (${\Bbb{G}},{\Bbb{J}},\Omega$) on 
the total space of the tangent bundle TN, given a positive definite metric g on the 2-manifold N. The associated
metric turned out to be of neutral signature and scalar-flat. The K\"ahler structure on the space ${\Bbb{L}}^{3}$ 
of oriented affine lines in ${\Bbb{E}}^3$, which can be identified with the tangent bundle to the 2-sphere, was
also investigated. 

In this paper we consider this K\"ahler structure on ${\Bbb{L}}^{3}_{1,-}=$TH$^2$. This K\"ahler structure 
has many properties in common with the ${\Bbb{E}}^3$ case and for this reason we state the shared features together in 
what follows. In particular, we prove that:

\vspace{0.1in}

\begin{IMThm}
The identity component of the isometry group of the K\"ahler metric on TS$^2$ (TH$^2$) is isomorphic to the identity 
component of the Euclidean (Lorentzian) isometry group.
\end{IMThm}

\vspace{0.1in}

This property has been proven to almost uniquely fix the K\"ahler structure on TS$^2$ \cite{salvai} and the same holds
for the K\"ahler structure on TH$^2$. We also investigate 
the geodesics of this metric:

\vspace{0.1in}

\begin{IMThm}
The geodesics of the K\"ahler metric on TS$^2$ (TH$^2$) are generated by the 1-parameter subgroups (transvections)
of the Euclidean (Lorentzian) isometry group.
\end{IMThm}

\vspace{0.1in}

Finally, we consider (space-like) Weingarten surfaces in ${\Bbb{E}}^3$ (${\Bbb{E}}^3_1$). These are surfaces
for which there exists a functional relationship between the eigenvalues of the 2nd fundamental form, and we
find the following characterisation for such surfaces:

\vspace{0.1in}

\begin{IMThm}
Let S be a $C^2$-smooth (space-like) surface  in ${\Bbb{E}}^3$ (${\Bbb{E}}^3_1$) and $\Sigma$ be the oriented normal 
congruence, considered as a surface in ${\Bbb{L}}^{3}$ $({\Bbb{L}}^{3}_{1,-})$. Then  S is Weingarten iff the Lorentzian metric 
induced by ${\Bbb{G}}$ on $\Sigma$ is scalar flat. 
\end{IMThm}

\vspace{0.1in}

In particular, S is minimal (maximal) in ${\Bbb{E}}^3$ (${\Bbb{E}}^3_1$) iff a certain holomorphic condition is 
satisfied, which we relate to the well-known Weierstrass representation \cite{hitch} \cite{kob}.

To prove the above results, we establish some general properties of the neutral K\"ahler metric in the next section.
We also consider the isometries, geodesics and 2-dimensional submanifolds of TN endowed with this metric.

Building on these results, Section 3 specialises to the metrics on TS$^2$ and TH$^2$ and contains the proofs of the 
Main Theorems on isometries, geodesics and 2-dimensional submanifolds. The final part of this section discusses the relationship with the Weierstrass representation.

\vspace{0.1in}

\section{Neutral K\"{a}hler Metric on TN}

Let (${\Bbb{M}},{\Bbb{G}},{\Bbb{J}},\Omega$) be a K\"ahler surface. That is, ${\Bbb{M}}$ is a real 4-manifold 
endowed with the following structures. First, there is the metric ${\Bbb{G}}$, which we do not insist be positive definite
- it may also have neutral signature ($++--$). In order to deal with both cases simultaneously we assume that the metric 
can be diagonalised pointwise to ($1,1,\epsilon,\epsilon$), for $\epsilon=\pm1$.

In addition, we have a complex structure ${\Bbb{J}}$, which is a mapping
${\Bbb{J}}:$T$_p{\Bbb{M}}\rightarrow $T$_p{\Bbb{M}}$ at each $p\in{\Bbb{M}}$, which satisfies ${\Bbb{J}}^2=-{\mbox{Id}}$. 
Finally, there is a symplectic form $\Omega$, which is a closed non-degenerate 2-form. These structures are required  
to satisfy the compatibility conditions:
\[
{\Bbb{G}}({\Bbb{J}}\cdot,{\Bbb{J}}\cdot)={\Bbb{G}}(\cdot,\cdot) \qquad\qquad
{\Bbb{G}}(\cdot,\cdot)=\Omega({\Bbb{J}}\cdot,\cdot).
\]

The following result highlights the difference between the case where ${\Bbb{G}}$ is positive definite and where 
it is neutral, and will be of use later.

\begin{Thm}\label{t:wirt}
Let $p\in{\Bbb{M}}$ and $v_1,v_2\in T_p{\Bbb{M}}$ span a plane. Then
\[
\Omega(v_1,v_2)^2+\epsilon\varsigma^2(v_1,v_2)={\mbox{det }}{\Bbb{G}}(v_i,v_j),
\]
where $\varsigma^2(v_1,v_2)\geq0$ with equality iff $\{v_1,v_2\}$ spans a complex plane.
\end{Thm}

\begin{pf}
Let us choose a unitary basis of frames $\{e_1,e_2\}$ with 
\[
{\Bbb{J}}(e_1)=ie_1 \qquad{\Bbb{J}}(e_2)=ie_2 \qquad\qquad {\Bbb{G}}(e_1,\bar{e}_1)=\epsilon  
                                 \qquad {\Bbb{G}}(e_2,\bar{e}_2)=1.
\]
In terms of the dual co-frame $\{\theta^1,\theta^2\}$ we can write
\[
{\Bbb{G}}=2\left(\epsilon\theta^1\bar{\theta}^1+\theta^2\bar{\theta}^2\right)
\qquad\qquad \Omega=i\left(\epsilon\theta^1\wedge\bar{\theta}^1+\theta^2\wedge\bar{\theta}^2\right).
\]
Thus
\[
\Omega\wedge\Omega=-2\epsilon\theta^1\wedge\bar{\theta}^1\wedge\theta^2\wedge\bar{\theta}^2.
\]
For $v_1,v_2\in T_p{\Bbb{M}}$, we evaluate this 4-form:
\begin{equation}\label{e:sympsq}
\Omega\wedge\Omega(v_1,v_2,{\Bbb{J}}(v_1),{\Bbb{J}}(v_2))
   =-2\epsilon\theta^1\wedge\bar{\theta}^1\wedge\theta^2\wedge\bar{\theta}^2(v_1,v_2,{\Bbb{J}}(v_1),{\Bbb{J}}(v_2)).
\end{equation}
Now, define the 2-form $\varsigma^2$ by
\[
\varsigma^2(v_1,v_2)={\textstyle{\frac{1}{4}}}\theta^1\wedge\bar{\theta}^1\wedge\theta^2\wedge\bar{\theta}^2(v_1,v_2,{\Bbb{J}}(v_1),{\Bbb{J}}(v_2)).
\]
If we decompose $v_1$ and $v_2$ on the unitary frame:
\[
v_j=a_je_1+\bar{a}_j\bar{e}_1+b_je_2+\bar{b}_j\bar{e}_2 \qquad\qquad a_j,b_j\in{\Bbb{C}}\qquad j=1,2,
\]
we compute that
\[
\varsigma^2(v_1,v_2)=-{\textstyle{\frac{1}{4}}}\left|
   \begin{matrix}
a_1 & \bar{a}_1 & b_1 & \bar{b}_1 \\
a_2 & \bar{a}_2 & b_2 & \bar{b}_2 \\
a_1 & -\bar{a}_1 & b_1 & -\bar{b}_1 \\
a_2 & -\bar{a}_2 & b_2 & -\bar{b}_2
\end{matrix}\right|=|a_1b_2-b_1a_2|^2\geq 0.
\]
It is clear that if the plane spanned by $v_1,v_2$ is holomorphic then $\varsigma^2$ vanishes. Conversely, suppose that
$\varsigma^2$ vanishes on some $v_1,v_2$ spanning a plane. Then for some $\mu\in{\Bbb{C}}$
\[
v_1=a_1e_1+b_1e_2+\bar{a}_2\bar{e}_1+\bar{b}_2\bar{e}_2
\qquad\qquad 
v_2=\mu(a_1e_1+b_1e_2)+\bar{\mu}(\bar{a}_2\bar{e}_1+\bar{b}_2\bar{e}_2) ,
\]
where $\mu$ is not real, since $v_1$ and $v_2$ are linearly independent. Thus we compute that
\[
{\Bbb{J}}(v_1)=i(a_1e_1+b_1e_2)-i(\bar{a}_2\bar{e}_1+\bar{b}_2\bar{e}_2) 
        =\frac{\mu+\bar{\mu}}{\bar{\mu}-\mu}iv_1-\frac{2}{\bar{\mu}-\mu}iv_2,
\]
and so the plane spanned by $v_1,v_2$ is holomorphic. 

Returning now to equation (\ref{e:sympsq}) we have
\begin{align}
\Omega(v_1,v_2)\Omega({\Bbb{J}}(v_1),{\Bbb{J}}(v_2))-&\Omega(v_2,{\Bbb{J}}(v_1))\Omega({\Bbb{J}}(v_2),v_1)\nonumber\\
   &-\Omega({\Bbb{J}}(v_2),v_2)\Omega(v_1,{\Bbb{J}}(v_1)) =-\epsilon\varsigma^2(v_1,v_2),\nonumber
\end{align}
or, 
\[
\Omega(v_1,v_2)^2+\epsilon\varsigma^2(v_1,v_2)={\Bbb{G}}(v_1,v_1){\Bbb{G}}(v_2,v_2)-{\Bbb{G}}(v_2,v_1){\Bbb{G}}(v_2,v_1),
\]
which yields the claimed result.

\end{pf}

We turn now to our construction of a neutral K\"ahler structure on TN - further details can be found in \cite{gak4}. 
Given a Riemannian 2-manifold (N,g,j) we construct
a canonical K\"ahler structure (${\Bbb{J}}$,$\Omega$,${\Bbb{G}}$) on the tangent bundle TN as follows. The
Levi-Civita connection associated with g splits the tangent bundle TTN$\cong$TN$\oplus$TN and the complex structure is
defined to be ${\Bbb{J}}=j\oplus j$. 

To define the symplectic form, consider the metric g as a mapping from TN to T$^*$N and pull back the canonical 
symplectic 2-form $\Omega^*$ on T$^*$N to a symplectic 2-form $\Omega$ on TN. 

Finally, the metric is defined as above by
${\Bbb{G}}(\cdot,\cdot)=\Omega({\Bbb{J}}\cdot,\cdot)$. The triple (${\Bbb{J}}$, $\Omega$, ${\Bbb{G}}$) 
determine a K\"ahler structure on TN. 

\begin{Prop}\cite{gak4}
Let (TN,${\Bbb{J}}$,$\Omega$,${\Bbb{G}}$) be the K\"ahler manifold, as above. Then the metric ${\Bbb{G}}$ has
neutral signature ($++--$) and is scalar-flat. Moreover, ${\Bbb{G}}$ is K\"ahler-Einstein iff g is flat, and
${\Bbb{G}}$ is conformally flat iff g is of constant curvature.
\end{Prop}

Choose conformal coordinates $\xi$ on N so that
$ds^2=e^{2u}d\xi d\bar{\xi}$, and corresponding coordinates ($\xi$,$\eta$) on TN  by
identifying
\[
(\xi,\eta) \leftrightarrow \eta\frac{\partial}{\partial \xi}+\bar{\eta}\frac{\partial}{\partial \bar{\xi}}
                   \in \mbox{T}_\xi N.
\] 

In such a coordinate system, the symplectic 2-form is

\[
\Omega=2{\Bbb{R}}\mbox{e}\left(e^{2u}d\eta\wedge d\bar{\xi}+\eta\partial(e^{2u})d\xi\wedge d\bar{\xi}\right),
\]
and the K\"ahler metric ${\Bbb{G}}$ is 
\[
{\Bbb{G}}=2{\Bbb{I}}\mbox{m}\left(e^{2u}d\bar{\eta} d\xi-\eta\partial(e^{2u})d\xi d\bar{\xi}\right).
\]
Here we have introduced the notation $\partial$ for differentiation with respect to $\xi$ - notation that we
will use throughout this paper.

\subsection{Isometries}

\vspace{0.1in}

\begin{Thm}\label{t:isom}
Let $Iso(TN,{\Bbb{G}})$ be the vector space of Killing vectors of ($TN,{\Bbb{G}}$) and $Iso(N,g)$ be the space of 
Killing vectors of ($N,g$), where the metric g is assumed to be complete. If g is non-flat, then
\[
Iso(TN,{\Bbb{G}})\cong Iso(N,g)\oplus\left({\cal{H}}ol(TN,N)\cap{\cal{L}}ag(TN,N)\right),
\]
where ${\cal{H}}ol(TN,N)$ and ${\cal{L}}ag(TN,N)$ are the spaces of holomorphic and Lagrangian sections 
of the canonical bundle $TN\rightarrow N$, respectively. In addition, dim ($Iso(N,g)$) =
dim (${\cal{H}}ol(TN,N)\cap{\cal{L}}ag(TN,N)$).

If g is flat, then
\[
Iso(TN,{\Bbb{G}})\cong Hty(N,g)\oplus\left({\cal{H}}ol(TN,N)\cap{\cal{L}}ag(TN,N)\right)\oplus {\cal{V}},
\]
where $Hty(N,g)$ is the space of homotheties of g and ${\cal{V}}$ is a certain 3-dimensional vector space.
\end{Thm}

\begin{pf}
A vector field ${\Bbb{V}}^i$ is a Killing vector for ${\Bbb{G}}$ if and only if
\[
{\Bbb{V}}^i\partial_i{\Bbb{G}}_{jk}+{\Bbb{G}}_{ki}\partial_j{\Bbb{V}}^i+{\Bbb{G}}_{ji}\partial_k{\Bbb{V}}^i=0.
\]

To find the solution to these equations we proceed as in Proposition 9 of \cite{gak4}. The equations with ($j,k$) 
equal to 
($\eta,\eta$), ($\eta,\bar{\eta}$) and ($\xi,\eta$) read: 
\[
\partial_\eta {\Bbb{V}}^{\bar{\xi}}=0
\qquad\qquad
\partial_\eta {\Bbb{V}}^{\xi}-\partial_{\bar{\eta}} {\Bbb{V}}^{\bar{\xi}}=0
\qquad\qquad
\partial_\xi {\Bbb{V}}^{\bar{\xi}}-\partial_{\eta} {\Bbb{V}}^{\bar{\eta}}=0,
\]
which have general solution
\begin{equation}\label{e:k}
{\Bbb{V}}^\xi=a_0(\xi,\bar{\xi}) +a_1(\xi,\bar{\xi})\eta
\qquad
{\Bbb{V}}^\eta=b_0(\xi,\bar{\xi},\eta)+\left(\bar{\partial} a_0 
       +  \eta\bar{\partial} a_1\right)\bar{\eta},
\end{equation}
for  complex functions $a_0(\xi,\bar{\xi})$ and $b_0(\xi,\bar{\xi},\eta)$, and  real function
$a_1(\xi,\bar{\xi})$. 

The $\eta$ derivative of the ($\xi,\bar{\eta}$) equation forces $a_1$ to be constant, while the remainder of the
equation says:
\[
\partial_{\bar{\eta}}\bar{b}_0+4a_1\bar{\partial}u\;\bar{\eta}+2a_0\partial u+2\bar{a}_0\bar{\partial}u
          +\partial_\xi a_0=0,
\]
which we integrate to
\[
b_0=b_1-(\partial_{\bar{\xi}} \bar{a}_0+2a_0\partial u+2\bar{a}_0\bar{\partial}u)\eta-2a_1\partial u\;\eta^2,
\]
for complex function $b_1=b_1(\xi,\bar{\xi})$.

With this simplification, we find that the $\bar{\eta}^2$-term says $a_1.\partial\bar{\partial} u=0$ and $b_1=b_1(\xi)$.
The Gauss curvature of the metric g is $\kappa=-4e^{-2u}\partial\bar{\partial}u$ and so we must now deal with two 
separate cases: $\kappa=0$ and $\kappa\neq0$. 

\vspace{0.1in}
\noindent{\bf Case 1}: $\kappa\neq0$

We have then that $a_1=0$, while the ($\xi$,$\xi$)-equation yields
\[
\partial\partial\bar{a}_0-2\partial u \partial \bar{a}_0=0
\]
\[
\partial\partial a_0+2\partial u \partial a_0-2a_0\partial\partial u+2\bar{a}_0\partial\bar{\partial}u=0.
\]
The first of these can be integrated to
\[
\bar{\partial}a_0=b_2e^{2u},
\]
where $b_2=b_2(\xi)$ is a complex-valued function. In fact, we must have $b_2=0$, since differentiating the  
($\xi$,$\bar{\xi}$)-equation with respect to $\bar{\xi}$ and $\bar{\eta}$ yields $b_2.\partial\bar{\partial} u=0$,
and $\partial\bar{\partial}u\neq0$. 
Thus $a_0$ is a holomorphic function and the second equation above can be integrated to
\begin{equation}\label{e:hom1}
\partial(a_0e^{2u})+\bar{\partial}(\bar{a}_0e^{2u})=C,
\end{equation}
for some real constant $C$. This equation can be understood as follows. Consider a vector field on the Riemannian 
2-manifold (N,g$=e^{2u}d\xi d\bar{\xi}$) given by
\[
X=a_0(\xi,\bar{\xi})\frac{\partial}{\partial \xi}+\bar{a}_0(\xi,\bar{\xi})\frac{\partial}{\partial \bar{\xi}}.
\]
Then $X$ is a homothety iff $a_0$ is holomorphic and satisfies equation (\ref{e:hom1}). Since the only complete
metrics admitting (non-Killing) homotheties are flat \cite{KN}, we conclude that $C=0$ and the vector field $X$ is 
in fact an infinitesimal isometry of (N,g).

The final equation to be solved is for the holomorphic function $b_1$:
\[
\partial(b_1e^{2u})=\bar{\partial}(\bar{b}_1e^{2u}).
\]
This is equivalent to the holomorphic section $\xi\rightarrow(\xi,\eta=b_1(\xi))$ being Lagrangian. The final 
form of the Killing vector is:
\[
{\Bbb{V}}^\xi=a_0\qquad\qquad
{\Bbb{V}}^\eta=b_1-(\partial_{\bar{\xi}} \bar{a}_0+2a_0\partial u+2\bar{a}_0\bar{\partial}u)\eta,
\]
where the holomorphic functions $a_0$ and $b_1$ correspond to Killing vectors of g and holomorphic
Lagrangian sections of TN (respectively). 

Note that we have an isomorphism $Iso(N,g)\cong{\cal{H}}ol(TN,N)\cap{\cal{L}}ag(TN,N)$ given by $a_0\rightarrow b_1=ia_0$.

\vspace{0.1in}
\noindent{\bf Case 2}: $\kappa=0$

We can take holomorphic coordinates $\xi$ on N such that $u=0$, and then the ($\xi$,$\xi$)-equation takes the form:
\[
\partial\partial a_0=0 \qquad\qquad \partial\partial \bar{a}_0=0 ,
\]
which has solution $a_0=b_4+b_5\xi+(b_6+b_7\xi)\bar{\xi}$ for complex constants $b_4,b_5,b_6,b_7$. The $\bar{\eta}$ term
of the ($\xi$,$\bar{\xi}$)-equation forces $b_7$ to vanish. The remaining equation to be solved is 
for the holomorphic function $b_1$:
\[
\partial(b_1)=\bar{\partial}(\bar{b}_1),
\]
with solution $b_1=b_2+b_3\xi$ for $b_2\in{\Bbb{C}}$ and $b_3\in{\Bbb{R}}$.
The resulting Killing vectors are
\[
{\Bbb{V}}^\xi=b_4+b_5\xi+b_6\bar{\xi}+a_1\eta \qquad\qquad
{\Bbb{V}}^\eta=b_2+b_3\xi-\bar{b}_5\eta+b_6\bar{\eta}.
\]
Thus the 10-dimensional space of Killing vectors consists of: the homotheties of g ($b_4,b_5\in{\Bbb{C}}$), the 
Lagrangian holomorphic sections of TN ($b_2\in{\Bbb{C}},b_3\in{\Bbb{R}}$) and the vector space ${\cal{V}}$ 
generated by $b_6\in{\Bbb{C}}$ and $a_1\in{\Bbb{R}}$.
\end{pf}

\subsection{Geodesics}

\vspace{0.1in}

We now consider geodesics of the neutral K\"ahler metric ${\Bbb{G}}$ on TN. 

\begin{Thm}
The linear subspaces of the fibres of the bundle TN$\rightarrow$N are null geodesics. These are the only
geodesics that lie in the fibres. 

The geodesics that do not lie in the fibres project under the bundle map to geodesics on N.
\end{Thm}
\begin{pf}
Consider the parameterized curve $s\rightarrow (\xi(s),\eta(s))$. The affinely parameterized geodesic equations 
for ${\Bbb{G}}$ are:
\[
\ddot{\xi}+2\partial u \dot{\xi}^2=0
\]
\[
\ddot{\eta}+4\partial u \dot{\xi}\dot{\eta}+2(\eta\partial\partial u -\bar{\eta}\bar{\partial}\partial u)\dot{\xi}^2 =0.
\]
If $\dot{\xi}(0)=0$, then the first equation implies that $\dot{\xi}(s)=0$ for all $s$ i.e. the geodesic lies in the 
fibre. The second equation can then be integrated to $\eta(s)=\eta(0)+\dot{\eta}(0)s$, and these geodesics
turn out to be null. This proves the first part of the theorem. 

On the other hand, if $\dot{\xi}(0)\neq0$, the first equation is just the geodesic equation for $g=e^{2u}d\xi d\bar{\xi}$ 
on N. This proves the second statement.
\end{pf}

In order to completely integrate the geodesic equations we need more information on g. In the next section we solve
the rest of these equations for g of constant curvature $\pm1$.

\subsection{Surfaces in TN}

\vspace{0.1in}

We now consider the geometry induced on an immersed surface $f:\Sigma\rightarrow$TN by the K\"ahler structure.
Let ($\nu,\bar{\nu}$) be a local coordinate system on $\Sigma$: 
$f(\nu,\bar{\nu})=(\xi(\nu,\bar{\nu}),\eta(\nu,\bar{\nu}))$. Then:

\begin{Prop}\label{p:lag}
\[
\Omega\left(f_*\frac{\partial}{\partial \nu},f_*\frac{\partial}{\partial \bar{\nu}}\right)=2e^{2u}{\Bbb{I}}m\left[
   \left(\partial_\nu\eta\partial_{\bar{\nu}}\bar{\xi}+\partial_\nu\bar{\eta}\partial_{\bar{\nu}}\xi
     \right)
     +2\eta\partial u\left(\partial_\nu\xi\partial_{\bar{\nu}}\bar{\xi}-\partial_{\bar{\nu}}\xi\partial_\nu\bar{\xi}\right)\right]
\]
\[
\varsigma^2\left(f_*\frac{\partial}{\partial \nu},f_*\frac{\partial}{\partial \bar{\nu}}\right)=
  e^{4u}|\partial_\nu\xi\partial_{\bar{\nu}}\eta-\partial_\nu\eta\partial_{\bar{\nu}}\xi|^2.
\]
\end{Prop}
\begin{pf}
The first of these comes from pulling back the symplectic 2-form to $\Sigma$, while the second follows from
\[
\varsigma^2\left(f_*\frac{\partial}{\partial \nu},f_*\frac{\partial}{\partial \bar{\nu}}\right)=
-{\textstyle{\frac{1}{4}}}e^{4u}\left|
   \begin{matrix}
\partial_\nu\xi & \partial_\nu\bar{\xi}& \partial_\nu\eta & \partial_\nu\bar{\eta} \\
\partial_{\bar{\nu}}\xi & \partial_{\bar{\nu}}\bar{\xi}& \partial_{\bar{\nu}}\eta & \partial_{\bar{\nu}}\bar{\eta} \\
\partial_\nu\xi & -\partial_\nu\bar{\xi}& \partial_\nu\eta & -\partial_\nu\bar{\eta} \\
\partial_{\bar{\nu}}\xi & -\partial_{\bar{\nu}}\bar{\xi}& \partial_{\bar{\nu}}\eta & -\partial_{\bar{\nu}}\bar{\eta} 
\end{matrix}\right|=e^{4u}|\partial_\nu\eta\partial_{\bar{\nu}}\xi-\partial_\nu\xi\partial_{\bar{\nu}}\eta|^2.
\]
\end{pf}

The signature of the induced metric on a Lagrangian surface is given by:

\begin{Prop}
The metric induced on a Lagrangian surface by the neutral K\"ahler metric is either Lorentzian or degenerate. This last occurs
when the surface is both Lagrangian and holomorphic.
\end{Prop}
\begin{pf}
A plane spanned by $\{v_1,v_2\}$ is Lagrangian if $\Omega(v_1,v_2)=0$. By Theorem \ref{t:wirt}, on a Lagrangian plane 
${\mbox{det }}{\Bbb{G}}(v_i,v_j)=-\varsigma^2(v_1,v_2)\leq0$. Thus the determinant of the metric induced on a Lagrangian 
surface is zero or negative depending on whether or not the surface is holomorphic.
\end{pf}

In the positive definite case this cannot occur: a surface can never be both holomorphic and Lagrangian, and the 
induced metric can be neither degenerate nor Lorentzian. In fact, the positive definite case of the Theorem \ref{t:wirt} 
leads to Wirtinger's inequality:
\[
\Omega(v_1,v_2)^2\leq{\mbox{det }}{\Bbb{G}}(v_i,v_j),
\]
with equality iff $\{v_1,v_2\}$ spans a complex plane.

The local expression for the induced metric is:

\begin{Prop}\label{p:indmet}
The metric induced on $f(\Sigma)$ by ${\Bbb{G}}$ is:
\begin{align}
f^*{\Bbb{G}}=&2e^{2u}{\Bbb{I}}\mbox{m}\left[\left( \partial_\nu\bar{\eta}\partial_\nu\xi
      -\partial_\nu\eta\partial_\nu\bar{\xi}
       -2(\eta\partial u-\bar{\eta}\bar{\partial} u)\partial_\nu\xi\partial_\nu\bar{\xi}\right)d\nu^2\right.\nonumber\\
&\qquad\qquad-\left.\left(\partial_\nu\eta\partial_{\bar{\nu}}\bar{\xi}
       +\partial_{\bar{\nu}}\eta\partial_\nu\bar{\xi}
   +2\eta\partial u(\partial_\nu\xi\partial_{\bar{\nu}}\bar{\xi}+\partial_{\bar{\nu}}\xi\partial_\nu\bar{\xi})\right)d\nu d\bar{\nu}\right].\nonumber
\end{align}
\end{Prop}
\begin{pf}
This follows from the definition of ${\Bbb{G}}$ pulled back by $f$.
\end{pf}

\vspace{0.1in}

\section{Oriented Lines in ${\Bbb{E}}^3$ and ${\Bbb{E}}^3_1$}

In the cases where N=S$^2$ or N=H$^2$ endowed with a metric of constant Gauss curvature 
($e^{2u}=4(1\pm\xi\bar{\xi})^{-2}$), 
the above construction yields the
neutral K\"ahler metric on the space ${\Bbb{L}}^{3}$ of oriented affine lines or on the space ${\Bbb{L}}^{3}_{1,-}$ of 
future-pointing time-like lines in ${\Bbb{E}}^3$ or ${\Bbb{E}}^3_1$ (respectively). 
Because of the similarity between these two cases we will treat them together - from here on reference to an oriented 
line will refer to just that in ${\Bbb{E}}^3$, and to a future-pointing time-like line in ${\Bbb{E}}^3_1$. The 
former case was dealt with in detail in \cite{gak4}.

With this in mind, we define the map $\Phi$ which sends ${\Bbb{L}}^{3}\times{\Bbb{R}}$ 
(${\Bbb{L}}^{3}_{1,-}\times{\Bbb{R}}$) to
${\Bbb{E}}^3$ (${\Bbb{E}}^3_1$) as follows: $\Phi$ takes an oriented line $\gamma$
and a real number $r$ to that point in ${\Bbb{E}}^3$ (${\Bbb{E}}^3_1$) which lies on $\gamma$ and is an affine parameter
distance $r$ from the point on $\gamma$ closest to the origin.

\begin{Prop}
The map can be written as $\Phi((\xi,\eta),r)=(z,t)\in{\Bbb{C}}\oplus{\Bbb{R}}={\Bbb{E}}^3$ (${\Bbb{E}}^3_1$) 
where the local coordinate expressions are:
\begin{equation}\label{e:coord1}
z=\frac{2(\eta\mp\bar{\eta}\xi^2)+2\xi(1\pm\xi\bar{\xi})r}{(1\pm\xi\bar{\xi})^2}
\qquad\qquad
t=\frac{\mp 2(\eta\bar{\xi}+\bar{\eta}\xi)+(1-\xi^2\bar{\xi}^2)r}{(1\pm\xi\bar{\xi})^2},
\end{equation}
with inverse
\begin{equation}\label{e:coord2}
\eta={\textstyle{\frac{1}{2}}}(z-2t\xi\mp\bar{z}\xi^2) \qquad\qquad 
          r=\frac{\pm\bar{\xi}z\pm\xi\bar{z}+(1\mp\xi\bar{\xi})t}{1\pm\xi\bar{\xi}},
\end{equation}
where the upper (lower) sign refers to the Euclidean (Lorentzian) case.
\end{Prop}
\begin{pf}
The proof of the Euclidean case was given in \cite{gak1}, while that of the Lorentzian case follows along 
almost identical lines.
\end{pf}

\subsection{Isometries}

\vspace{0.1in}

\begin{Main}
The identity component of the isometry group of the metric on TS$^2$ (TH$^2$) is isomorphic to the identity component 
of the Euclidean (Lorentzian) isometry group.
\end{Main}
\begin{pf}
This was proven for TS$^2$ in \cite{gak4}, so we prove only the Lorentzian case here. The identity component 
of the (affine) isometries of ${\Bbb{E}}^3_1$ consists of rotations and translations:
${\mbox{ISO}}_0({\Bbb{E}}^3_1)\cong\mbox{SO}_0(2,1)\ltimes{\Bbb{E}}^3_1$. We first consider the action of
this group on ${\Bbb{L}}^{3}_{1,-}$. 

The 3-dimensional group of rotations about the origin preserve the hyperbolic plane and are given by the fractional
linear transformations
\[
\xi\rightarrow\xi'=\frac{\alpha\xi+\beta}{\bar{\beta}\xi+\bar{\alpha}}
\qquad\qquad {\mbox{for}} \qquad \alpha\bar{\alpha}-\beta\bar{\beta}=1.
\]
The group action on the space of future-pointing time-like lines is the derivative of this action:
\[
(\xi,\eta)\rightarrow(\xi',\eta')=\left(\frac{\alpha\xi+\beta}{\bar{\beta}\xi+\bar{\alpha}}
  ,\frac{1}{(\bar{\beta}\xi+\bar{\alpha})^2}\eta\right).
\]
The derivative at the identity of this action gives the infinitesimal generators:
\[
{\Bbb{V}}={\Bbb{R}}{\mbox{e}}\left[\left(\dot{\beta}+2\dot{\alpha}\xi-\dot{\bar{\beta}}\xi^2\right)
    \frac{\partial}{\partial\xi}-2\left(\dot{\bar{\alpha}}+\dot{\bar{\beta}}\xi\right)
        \eta\frac{\partial}{\partial\eta}\right],
\]
where $\dot{\bar{\alpha}}=-\dot{\alpha}$.

On the other hand, the infinitesimal translations act on ${\Bbb{L}}^{3}_{1,-}$ by ({\it cf.} equation (\ref{e:coord2})):
\[
{\Bbb{V}}={\Bbb{R}}{\mbox{e}}\left[{\textstyle{\frac{1}{2}}}(\dot{\gamma}-2\dot{\delta}\xi+\dot{\bar{\gamma}}\xi^2) 
      \frac{\partial}{\partial\eta}\right],
\]
where $\dot{\bar{\delta}}=\dot{\delta}$.

We now apply Theorem \ref{t:isom} to TH$^2$ and see that the infinitesimal isometries of ${\Bbb{G}}$ on 
${\Bbb{L}}^{3}_{1,-}$  are exactly the infinitesimal isometries of the Lorentzian metric on ${\Bbb{E}}^3_1$, with 
the identification:
\[
a_0=\dot{\beta}+2\dot{\alpha}\xi-\dot{\bar{\beta}}\xi^2
\qquad\qquad
b_1={\textstyle{\frac{1}{2}}}(\dot{\gamma}-2\dot{\delta}\xi+\dot{\bar{\gamma}}\xi^2).
\]
The result follows.
\end{pf}

In \cite{salvai} it is proven that above metric is the only K\"ahler metric on TS$^2$ that is invariant under the action 
induced on TS$^2$ (considered as the space of oriented lines in ${\Bbb{E}}^3$) by the Euclidean group 
(up to the addition of the pull-back of the round metric on S$^2$). A similar result holds for TH$^2$: the above metric 
is the only K\"ahler metric on TH$^2$ that is invariant under the action 
induced on TH$^2$ (considered as the space of future pointing time-like lines in ${\Bbb{E}}^3_1$) by the Lorentz group 
(up to the addition of the pull-back of the hyperbolic metric on H$^2$).

\vspace{0.1in}

\subsection{Geodesics}

\vspace{0.1in}

\begin{Main}
The geodesics of the metric ${\Bbb{G}}$ on TS$^2$ (TH$^2$) are generated by the 1-parameter subgroups (transvections)
of the Euclidean (Lorentzian) isometry group.
\end{Main}
\begin{pf}
In \cite{gak4} this was proven for the Euclidean case and so we concentrate on the Lorentzian case here. 
Let $c:[0,1]\rightarrow TH^2$ be a curve with tangent vector
\[
{\Bbb{X}}=\dot{\xi}\frac{\partial}{\partial \xi}
    +\dot{\eta}\frac{\partial}{\partial \eta}
   +\dot{\bar{\xi}}\frac{\partial}{\partial \bar{\xi}}
   +\dot{\bar{\eta}}\frac{\partial}{\partial \bar{\eta}}.
\]
We can set $\xi(0)=\eta(0)=0$, so that
the initial line $c(0)$ is the $x^3$-axis. We still
retain the freedom to rotate about, and translate along, the $x^3$-axis.

The geodesic equations (with arc-length or affine parameter $s$) are

\begin{equation}\label{e:geo1}
\ddot{\xi}+\frac{2\bar{\xi}}{1-\xi\bar{\xi}}\dot{\xi}^2=0
\end{equation}
\begin{equation}\label{e:geo2}
\ddot{\eta}+\frac{4\bar{\xi}}{1-\xi\bar{\xi}}\dot{\xi}\dot{\eta}+\frac{2(\bar{\xi}^2\eta-\bar{\eta})}{(1-\xi\bar{\xi})^2}\dot{\xi}^2=0,
\end{equation}
where a dot represents differentiation with respect to $s$, and we have
made use of the connection coefficients associated with ${\Bbb{G}}$. 

These have first integral
\begin{equation}\label{e:geo3}
\frac{2i}{(1-\xi\bar{\xi})^2}\left(
  \dot{\eta} \dot{\bar{\xi}}-\dot{\bar{\eta}} \dot{\xi}
   -\frac{2(\xi\bar{\eta}-\bar{\xi}\eta)}{1-\xi\bar{\xi}}\dot{\xi} \dot{\bar{\xi}}\right)=C_1,
\end{equation}
where $C_1\in{\Bbb{R}}$ vanishes iff the geodesic is null.

We have already considered the geodesics which lie in the fibre, and so we assume now that $\dot{\xi}(0)\neq 0$. 
Equation (\ref{e:geo1}) is the geodesic equation on the hyperbolic plane.  We can integrate this 
(with initial condition $\xi(0)=0$) to determine the evolution of $\xi$:
\begin{equation}\label{e:soln1}
\xi=\tanh(C_2s)e^{i\theta},
\end{equation}
for constants $C_2\in{\Bbb{R}}$ and $\theta\in[0,2\pi)$. Substituting this in
(\ref{e:geo3}) we find that
\begin{equation}\label{e:geo3a}
\eta e^{-i\theta}-\bar{\eta} e^{i\theta}
      =\frac{C_1s+C_3}{2iC_2\cosh^2(C_2s)},
\end{equation}
for some real constant $C_3$. Equation (\ref{e:geo2}) now simplifies
to:
\[
\frac{d^2}{ds^2}(\cosh^2(C_2s)\eta)
   -4C_2^2\cosh^2(C_2s)\eta =C_2(C_1s+C_3)ie^{i\theta},
\] 
with solution (noting equation (\ref{e:geo3a}))
\[
\eta=\frac{C_4\cosh(2C_2s)+C_5\sinh(2C_2s)-(C_1s+C_3)i}{4C_2\cosh^2(C_2s)}
    e^{i\theta},
\]
for real constants $C_4$ and $C_5$.
Finally, since $\eta(0)=0$ we have $C_3=C_4=0$, so that
\begin{equation}\label{e:soln2}
\eta=\frac{C_5\sinh(2C_2s)-C_1si}{4C_2\cosh^2(C_2s)}
    e^{i\theta}.
\end{equation}
Equations (\ref{e:soln1}) and (\ref{e:soln2}) are
the general solution to the geodesic equations in TH$^2$ when
the initial line is the $x^3$-axis. The four real
constants remaining, namely $C_1,C_2,C_5$ and $\theta$, determine the
initial direction of the geodesic in TH$^2$. 

The associated ruled surface in ${\Bbb{E}}^3$ is
a hyperbolic helicoid when $C_1\neq 0$ and a plane when $C_1=0$. To see this we
can put it in standard position as follows. By a rotation about the
$x^3$-axis we can fix $\theta=0$, while a translation along the $x^3$-axis 
allows us to put $C_5=0$. The ruled surface can be explicitly
determined using (\ref{e:coord1}) and the result
is
\[
x^1=t\sinh(2C_2s) \qquad x^2=-\frac{C_1s}{2C_2} \qquad x^3=t\cosh(2C_2s). 
\]
This is a hyperbolic helicoid for $C_1\neq 0$ and a plane for $C_1=0$, as claimed. These are
precisely the transvections of the Lorentzian isometry group.
\end{pf}

\subsection{Line congruences}

\vspace{0.1in}

A  {\it line congruence} is a 2-parameter family of oriented lines in ${\Bbb{E}}^3$ or ${\Bbb{E}}^3_1$ - that is,
a surface in ${\Bbb{L}}^{3}$ or ${\Bbb{L}}^{3}_{1,-}$. Locally, such a line congruence is given parametrically by a map
from ${\Bbb{C}}$ to ${\Bbb{C}}^2$: $\nu\rightarrow(\xi(\nu,\bar{\nu}),\eta(\nu,\bar{\nu}))$.

Given 
a line congruence, away from crossings, we can construction an adapted null frame: that is a trio of complex vector fields
$\{e_{(0)}, e_{(+)}, e_{(-)}\}$ such that $e_{(0)}=\overline{e_{(0)}}$, $e_{(+)}=\overline{e_{(-)}}$,
$\pm e_{(0)}\cdot e_{(0)}=e_{(+)}\cdot e_{(-)}=1$, and $e_{(0)}\cdot
e_{(+)}=0$, where the Euclidean or Lorentzian inner product (denoted by a dot) is extended bilinearly over 
${\Bbb{C}}$ and $e_{(0)}$ 
is aligned with the direction of the lines. Let $\{\theta^{(0)}, \theta^{(+)}, \theta^{(-)}\}$
be the dual coframes.

The complex {\it spin coefficients} $\Gamma_{m\;\;n}^{\;\;\;p}$ of the congruence are defined by
\[
e_{(p)}^j\nabla_je_{(m)}^{\;\;\;i}=\Gamma_{m\;\;p}^{\;\;n}e_{(n)}^i,
\]
where $\nabla$ is the flat covariant derivative and the indices
$m$, $n$, $p$ range over $0$, $+$, $-$. Breaking covariance, we introduce
the complex {\it optical scalars}: 
\[
\Gamma_{+\;\;-}^{\;\;\;0}=\rho \qquad\qquad \Gamma_{+\;\;+}^{\;\;\;0}=\sigma. 
\]

The following two Propositions establish the geometric significance of these complex functions.

\begin{Prop}
A line congruence is orthogonal to a (space-like) surface in ${\Bbb{E}}^3$ (${\Bbb{E}}^3_1$) iff $\rho$ is real.
\end{Prop}
\begin{pf}
The distribution of planes orthogonal to the line congruence are integrable iff $[e_+,e_-]\subset{\mbox{span}}\{e_+,e_-\}$
where $[\cdot,\cdot]$ is the Lie bracket of vector fields. Equivalently, integrability is $e_0\cdot[e_+,e_-]=0$.

Now since the connection is torsion-free we have that
\[
[e_m,e_n]=\left(\Gamma_{n\;\;m}^{\;\;\;p}-\Gamma_{m\;\;n}^{\;\;\;p}\right)e_p.
\]
Thus
\[
e_0\cdot[e_+,e_-]=\pm\Gamma_{-\;\;+}^{\;\;\;0}\mp\Gamma_{+\;\;-}^{\;\;\;0}=\pm(\bar{\rho}-\rho).
\]
The result follows.
\end{pf}

\begin{Prop}\label{p:2ndfund}
In the case where the line congruence is orthogonal to a surface S in ${\Bbb{E}}^3$ or ${\Bbb{E}}^3_1$, these
functions have the following expressions in terms of the eigenvalues $\lambda_1,\lambda_2$ of the 2nd fundamental form
of S:
\[
|\sigma|^2={\textstyle{\frac{1}{4}}}\left(\lambda_1-\lambda_2\right)^2
\qquad\qquad
\rho^2={\textstyle{\frac{1}{4}}}\left(\lambda_1+\lambda_2\right)^2.
\]
In addition, the argument of $\sigma$ determines the eigen-directions of the 2nd fundamental form. 
\end{Prop}

\begin{pf}
Let S be a C$^2$ (space-like) surface immersed in ${\Bbb{E}}^3$ (${\Bbb{E}}^3_1$) and let $N_i$ be
the unit normal. The second fundamental form is a symmetric two tensor on S
defined by
\[
h_{ij}=P_i^kP_j^l\nabla_kN_l,
\]
where $\nabla$ is the flat connection on ${\Bbb{E}}^3$ or ${\Bbb{E}}^3_1$ and $P_i^j$ is orthogonal projection onto
the tangent space of S. Such a symmetric two tensor has two real eigenvalues $\lambda_1$ and
$\lambda_2$ at each point of S. These are called the {\it principal  curvatures} of S and the
associated eigen-directions are called {\it principal directions} for the surface. Let $\{e_1,e_2\}$ be an orthonormal
eigen-basis for $H_{ij}$, and then the null frame is $e_{+}={\textstyle{\frac{1}{\sqrt{2}}}}(e_1+ie_2)e^{i\alpha}$ for
some angle $\alpha$.
Computing the null frame components of $H_{ij}$
\[
H_{++}=<\nabla_{e_+}\theta^0,e_+>=\Gamma_{+\;\;+}^{\;\;\;0}=\sigma\qquad
H_{+-}=<\nabla_{e_+}\theta^0,e_->=\Gamma_{+\;\;-}^{\;\;\;0}=\bar{\rho}=\rho.
\]
In terms of the real basis
\[
\sigma=H_{++}={\textstyle{\frac{1}{2}}}\left(H_{11}-H_{22}\right)e^{2i\alpha} 
   ={\textstyle{\frac{1}{2}}}\left(\lambda_1-\lambda_2\right)e^{2i\alpha}
\]
\[
\rho=H_{+-}={\textstyle{\frac{1}{2}}}\left(H_{11}+H_{22}\right) 
   ={\textstyle{\frac{1}{2}}}\left(\lambda_1+\lambda_2\right).
\] 
Finally, the argument of the shear is equal to $2\alpha$ which measures the angle of rotation between the frame and the
eigen-directions.
\end{pf}

These functions are given for a parametric line congruence by:

\begin{Thm}\label{t:spinco}

Let $\nu\rightarrow(\xi(\nu,\bar{\nu}),\eta(\nu,\bar{\nu}))$ be a parametric line congruence and define
\[
\partial^+\eta=\partial\eta+r\partial\xi\mp\frac{2\bar{\xi}\eta}{1\pm\xi\bar{\xi}}\partial\xi
\qquad\qquad
\partial^-\eta=\bar{\partial}\eta+r\bar{\partial}\xi\mp\frac{2\bar{\xi}\eta}{1\pm\xi\bar{\xi}}\bar{\partial}\xi,
\]
where $\partial$ is differentiation with respect to $\nu$, and the upper sign refers to ${\Bbb{E}}^3$ and the 
lower to ${\Bbb{E}}^3_1$.

Then
\[
\rho=\frac{ \partial^+\eta\bar{\partial}\;\bar{\xi} -\partial^-\eta\partial\bar{\xi}}
{\partial^-\eta\overline{\partial^-\eta}-\partial^+\eta\overline{\partial^+\eta}}
\qquad\qquad
\sigma=\frac{\overline{\partial^+\eta}\partial\bar{\xi} -\overline{\partial^-\eta}\;\bar{\partial}\;\bar{\xi}}
{\partial^-\eta\overline{\partial^-\eta}-\partial^+\eta\overline{\partial^+\eta}}.
\]
\end{Thm}

\begin{pf}
The following is an adapted null frame:
\[
e_0=D\Phi\left(\frac{\partial}{\partial r}\right)
\qquad\qquad
e_+= \alpha D\Phi\left(\frac{\partial}{\partial \nu}\right)
     +\beta D\Phi\left(\frac{\partial}{\partial \bar{\nu}}\right)
     +\Omega D\Phi\left(\frac{\partial}{\partial r}\right),
\]
where
\[
\Omega=\frac{\sqrt{2}\left[\overline{\partial^-\eta}(\eta\bar{\partial}\;\bar{\xi}+\bar{\eta}\;\bar{\partial}\xi)-\overline{\partial^+\eta}(\eta\partial\bar{\xi}+\bar{\eta}\partial\xi)\right]}
{(1\pm\bar{\xi}\xi) (\partial^-\eta\overline{\partial^-\eta}-\partial^+\eta\overline{\partial^+\eta}) },
\]
and
\[
\alpha=\mp\frac{\overline{\partial^+\eta}(1\pm\bar{\xi}\xi)}
{\sqrt{2}(\partial^-\eta\overline{\partial^-\eta}-\partial^+\eta\overline{\partial^+\eta}) }
\qquad
\beta=\pm\frac{\overline{\partial^-\eta}(1\pm\bar{\xi}\xi)}
{\sqrt{2}(\partial^-\eta\overline{\partial^-\eta}-\partial^+\eta\overline{\partial^+\eta}) }.
\]
This can be checked by computing the derivative of $\Phi$ and seeing that 
$\pm e_{(0)}\cdot e_{(0)}=e_{(+)}\cdot e_{(-)}=1$, and $e_{(0)}\cdot e_{(+)}=0$.

Now we introduce the complex vectors:
\[
Z_+=D\Phi\left(\frac{\partial}{\partial\nu}\right)
-e_0\cdot D\Phi\left(\frac{\partial}{\partial\nu}\right)D\Phi\left(\frac{\partial}{\partial r}\right)
\]
\[
Z_-=D\Phi\left(\frac{\partial}{\partial\bar{\nu}}\right)
-e_0\cdot D\Phi\left(\frac{\partial}{\partial\bar{\nu}}\right)D\Phi\left(\frac{\partial}{\partial r}\right).
\]
These have inner product
\[
 Z_+\cdot Z_+=\frac{4\partial^+\eta\overline{\partial^-\eta}}{(1\pm\xi\overline{\xi})^2}
\qquad
Z_+\cdot Z_-=\frac{2(\partial^-\eta\overline{\partial^-\eta}+\partial^+\eta\overline{\partial^+\eta})}
 {(1\pm\xi\overline{\xi})^2}.
\]
Then
\[
\bar{\rho} = [e_0,e_+]\cdot e_-= \left(\frac{\partial \alpha}{\partial r}Z_++\frac{\partial \beta}{\partial r}Z_-\right)
  \cdot\left(\bar{\alpha}Z_-+\bar{\beta}Z_+\right),
\]
and
\[
\sigma = [e_0,e_+]\cdot e_+= \left(\frac{\partial \alpha}{\partial r}Z_++\frac{\partial \beta}{\partial r}Z_-\right)
   \cdot\left(\alpha Z_++\beta Z_-\right).
\]
Now use the expressions for $\alpha$ and $\beta$ and the equation above to get the result.
\end{pf}

A surface S in ${\Bbb{E}}^3$ (${\Bbb{E}}^3_1$) is said to be {\it Weingarten} if the eigenvalues of the 
second fundamental form of S are functionally dependent.

\vspace{0.1in}

\begin{Main}
Let S be a $C^2$-smooth (space-like) surface  in ${\Bbb{E}}^3$ (${\Bbb{E}}^3_1$) and $\Sigma$ be the oriented normal 
congruence, considered as a surface in ${\Bbb{L}}^{3}$ $({\Bbb{L}}^{3}_{1,-})$. Then  S is Weingarten iff the Lorentzian metric 
induced by ${\Bbb{G}}$ on $\Sigma$ is scalar flat. 
\end{Main}
\begin{pf}
Let S be a $C^2$ (space-like) surface in ${\Bbb{E}}^3$ (${\Bbb{E}}^3_1$) and $f:\Sigma\rightarrow{\Bbb{L}}^{3}$ 
(${\Bbb{L}}^{3}_{1,-}$) the oriented normal congruence. Let $\nu\rightarrow(\xi(\nu,\bar{\nu}),\eta(\nu,\bar{\nu}))$ 
be a local parameterization of this line congruence. 

The lagrangian condition is an
integrability condition (locally) for a real function, which can be identified with distance $r=r(\nu,\bar{\nu})$
of the point on the surface S from the closest point to the origin along the oriented normal line. In fact:
\begin{equation}\label{e:supfunc}
\partial_{\bar{\nu}}r=\pm\frac{2\eta\partial_{\bar{\nu}}\bar{\xi}
                            +2\bar{\eta}\partial_{\bar{\nu}}\xi}{(1\pm\xi\bar{\xi})^2}.
\end{equation}

First, let us assume that the Weingarten surface S is flat: $\lambda_1\lambda_2=0$. In this case \cite{gak2}
$\partial_\nu\xi\partial_{\bar{\nu}}\bar{\xi}-\partial_{\bar{\nu}}\xi\partial_\nu\bar{\xi}=0$ and we can 
parameterize the surface by the fibre coordinate $\nu=\eta$. From Proposition \ref{p:lag}
\[
f^*\Omega=e^{2u}\left(\partial_\nu\xi-\partial_{\bar{\nu}}\bar{\xi}\right)d\nu\wedge d\bar{\nu},
\]
and since $\Sigma$ is Lagrangian $\partial_\nu\xi=\partial_{\bar{\nu}}\bar{\xi}$.
Then, the induced metric is ({\it cf}. Proposition \ref{p:indmet}) 
\[
f^*{\Bbb{G}}=\frac{8}{(1\pm\xi\bar{\xi})^2}{\Bbb{I}}{\mbox{m}}\left[\left(-\partial_\nu\bar{\xi}\pm\frac{2(\nu\bar{\xi}-\bar{\nu}\xi)}{1\pm\xi\bar{\xi}}\partial_\nu\xi\partial_\nu\bar{\xi}\right)d\nu^2
   \pm\frac{4\nu\bar{\xi}}{1\pm\xi\bar{\xi}}(\partial_\nu\xi)^2d\nu d\bar{\nu}
\right].
\]
A direct computation then shows that the scalar curvature ${\Bbb{K}}$ of this metric vanishes, as claimed.

On the other hand, suppose the curvature of S is not zero, then the oriented normal congruence is the graph of a
section of the canonical bundle TN$\rightarrow$N. Such surfaces are characterised by the fact that
$\partial_\nu\xi\partial_{\bar{\nu}}\bar{\xi}-\partial_{\bar{\nu}}\xi\partial_\nu\bar{\xi}\neq0$. In such 
a case, it is natural to parameterize the surface by the coordinate on the base: $\nu=\xi$.

Let us write this section as $\xi\rightarrow(\xi,\eta=F(\xi,\bar{\xi}))$ and denote the complex 
slopes of the section of TN$\rightarrow$N by: 
\[
\sigma_0=-\partial\bar{F} \qquad\qquad \rho_0=e^{-2u}\partial (e^{2u}F).
\]
The functions $\rho_0$, $\sigma_0$ and $u$ satisfy the following differential relation:
\begin{equation}\label{e:cm2}
\bar{\partial}\rho_0=-e^{-2u}\partial(\bar{\sigma}_0e^{2u})-{\textstyle{\frac{1}{2}}}Fe^{2u}\kappa,
\end{equation}
where $\kappa$ is the Gauss curvature of  the metric g on N.
This follows from the fact that partial derivatives commute. The left hand side reads
\[
\bar{\partial}\rho_0=\bar{\partial}e^{-2u}\partial (e^{2u}F)=\bar{\partial}\partial F+2\partial u\bar{\partial}F
    +2F\bar{\partial}\partial u,
\]
while the right hand side is
\[
-e^{-2u}\partial(\bar{\sigma}_0e^{2u})-{\textstyle{\frac{1}{2}}}Fe^{2u}\kappa=
e^{-2u}\partial(\bar{\partial}Fe^{2u})+2F\bar{\partial}\partial u=
\partial\bar{\partial}F+2\partial u\bar{\partial}F+2F\bar{\partial}\partial u,
\]
where we have used the expression for the Gaussian curvature in conformal coordinates: 
$\kappa=-4e^{-2u}\bar{\partial}\partial u$.

It follows from Proposition \ref{p:lag} that:
\[
\Omega\left(f_*\frac{\partial}{\partial \xi},f_*\frac{\partial}{\partial \bar{\xi}}\right)=
2{\Bbb{I}}m\left[e^{2u}\partial F+F\partial e^{2u}\right]=2e^{2u}{\Bbb{I}}m(\rho_0).
\]
Since $\Sigma$ is a Lagrangian section, $\rho_0$ is real:
\[
\partial\left(\frac{F}{(1\pm\xi\bar{\xi})^2}\right)=\bar{\partial}\left(\frac{\bar{F}}{(1\pm\xi\bar{\xi})^2}\right),
\]
and equation (\ref{e:supfunc}) simplifies to the integrability condition of the above equation:
\begin{equation}\label{e:supfunc2}
\bar{\partial} r=\pm\frac{2F}{(1\pm\xi\bar{\xi})^2}.
\end{equation}
By Proposition \ref{p:indmet}, the metric induced on a Lagrangian graph is
\[
ds^2=2ie^{2u}\left(\sigma_0 d\xi^2-\bar{\sigma}_0 d\bar{\xi}^2\right).
\]
Computing the curvature of this metric we get
\[
{\Bbb{K}}=\frac{e^{-4u}}{4|\sigma_0|^4}{\Bbb{I}}m\left[ \partial(|\sigma_0|^{2})\partial(\bar{\sigma}_0e^{2u})
        +Fe^{4u}|\sigma_0|^2\partial\kappa\right],
\]
where there are no second order terms in $\sigma_0$ because of the following identity, obtained by differentiation of 
equation (\ref{e:cm2}):
\[
{\Bbb{I}}m\left[\partial\left(e^{-2u}\partial(\bar{\sigma}_0e^{2u})\right)+{\textstyle{\frac{1}{2}}}Fe^{2u}\partial\kappa
    \right]=0.
\]
Substitution of equation (\ref{e:cm2}) yields the scalar curvature of the graph of a Lagrangian section of 
the canonical bundle to be:
\[
{\Bbb{K}}=-\frac{e^{-2u}}{2|\sigma_0|^4}{\Bbb{I}}m\left[ \partial\left(|\sigma_0|^{2}\right)\bar{\partial}\rho_0
        +{\textstyle{\frac{1}{2}}}Fe^{2u}\partial\left(|\sigma_0|^{2}\kappa\right)\right].
\]

In our case $N=S^2(H^2)$ and $\kappa=\pm1$ and so, with the aid of equation (\ref{e:supfunc2}), 
the scalar curvature simplifies to  
\[
{\Bbb{K}}=-\frac{(1\pm\xi\bar{\xi})^2}{8|\sigma_0|^4}{\Bbb{I}}\mbox{m}
     \left[ \partial\left(|\sigma_0|^{2}\right)\bar{\partial}\left(r+\rho_0\right)\right].
\]
Finally the spin coefficients are related to the slopes of the graph by Theorem \ref{t:spinco} (with $\nu=\xi$):
\[
\sigma=\frac{\sigma_0}{(r+\rho_0)^2-|\sigma_0|^2}
\qquad\qquad
\rho=\frac{r+\rho_0}{(r+\rho_0)^2-|\sigma_0|^2}.
\]

For a Lagrangian line congruence it follows from Proposition \ref{p:2ndfund} that
\[
|\sigma_0|^2={\textstyle{\frac{1}{4}}}\left({\textstyle{\frac{1}{\lambda_1}}}-{\textstyle{\frac{1}{\lambda_2}}}\right)^2
\qquad\qquad
(r+\rho_0)^2={\textstyle{\frac{1}{4}}}\left({\textstyle{\frac{1}{\lambda_1}}}+{\textstyle{\frac{1}{\lambda_2}}}\right)^2,
\]
where $\lambda_1$ and $\lambda_2$ are the eigenvalues of the second fundamental form of the surface S 
in ${\Bbb{E}}^3$ (${\Bbb{E}}^3_1$) that is orthogonal to the lines. Then if S is Weingarten, we
have $d\lambda_1\wedge d\lambda_2=0$ and the scalar curvature of the line congruence in ${\Bbb{L}}^{3}$ $({\Bbb{L}}^{3}_{1,-})$ 
vanishes.

Conversely, let $\Sigma$ be a Lagrangian line congruence in TS$^2$ (TH$^2$) with zero scalar curvature. If $\Sigma$ is
not the graph of a section, then it is orthogonal to a flat surface S in ${\Bbb{E}}^3$ (${\Bbb{E}}^3_1$), 
which is certainly Weingarten. On the 
other hand, if $\Sigma$ is the graph of a section, then, by the above computation, the vanishing of the scalar curvature
implies  $d\lambda_1\wedge d\lambda_2=0$ and once again, the orthogonal surface S in ${\Bbb{E}}^3$ (${\Bbb{E}}^3_1$) 
is Weingarten.
\end{pf}
\vspace{0.1in}

\subsection{Minimal (Maximal) Surfaces in ${\Bbb{E}}^3$ (${\Bbb{E}}^3_1$)}
In the special case of minimal (maximal) surfaces the mean curvature vanishes, that is, $\rho=0$.

\begin{Thm}

A Lagrangian congruence $\Sigma\subset TN$, for $N=S^2$ (H$^2$) is orthogonal to a
minimal (maximal) surface without flat points iff the congruence is the graph of
a local section $\xi\mapsto(\xi,\eta=F(\xi,\bar{\xi}))$ with
\begin{equation}\label{e:mineq}
\bar{\partial}\left(\frac{\partial\bar{F}}{(1\pm\xi\bar{\xi})^2}\right)=0.
\end{equation}
\end{Thm}
\begin{pf}
Let S be a minimal (maximal) surface without flat points and $\Sigma$ be its normal line congruence which is given
by the graph of a section. 
Then, by Proposition \ref{p:2ndfund} $\rho=0$ on S and thus we must have $r+\rho_0=0$. Now the identity 
(\ref{e:cm2}) can be written
\begin{equation}\label{e:cm3}
\bar{\partial}(r+\rho_0)=-(1\pm\xi\bar{\xi})^2\partial\left(\frac{\bar{\sigma}_0}{(1\pm\xi\bar{\xi})^2}\right).
\end{equation}
Putting $r+\rho_0=0$ in this identity yields the result.

Conversely, suppose equation (\ref{e:mineq}) holds for a Lagrangian line congruence $\Sigma$ which is given by 
the graph of a local section. Then, by equation (\ref{e:cm3}) $r+\rho_0=C$ for some real constant $C$. 
As the orthogonal surfaces move along the line congruence in ${\Bbb{E}}^3$ (${\Bbb{E}}^3_1$) we find that
$r\rightarrow r+{\mbox{ constant}}$. Thus there exists a surface S for which $r+\rho_0=0$, and therefore
$\rho=0$, i.e. there is a minimal surface orthogonal to $\Sigma$.
\end{pf}

The previous Theorem has two immediate consequences:

\begin{Cor}
The normal congruence to a minimal surface is given (up to translation) by a local
 section $F$ of the bundle $\pi:TN\rightarrow N$ with
\[
F=\sum_{n=0}^\infty 2\lambda_n\xi^{n+3}-\bar{\lambda}_n\bar{\xi}^{n+1}
  \left({\scriptstyle{\pm(n+2)(n+3)}}+2{\scriptstyle{(n+1)(n+3)}}\xi\bar{\xi}\pm{\scriptstyle{(n+1)(n+2)}}\xi^2\bar{\xi}^2
\right),
\]
for complex constants $\lambda_n$. The potential function $r:\Sigma\rightarrow{\Bbb{R}}$ satisfying (\ref{e:supfunc2}) is:
\[
r=-2\sum_{n=0}^\infty \frac{(3+n\pm(1+n)\xi\bar{\xi})(\lambda\xi^{n+2}+\bar{\lambda}\bar{\xi}^{n+2})}{1\pm\xi\bar{\xi}},
\]
where the upper (lower) sign refers to the Euclidean (Lorentzian) case.
\end{Cor}
\begin{pf}
Since the minimal surface condition is a holomorphic condition we can expand in a power series about a point:
\[
\frac{\partial\bar{F}}{(1\pm\xi\bar{\xi})^2}=\sum_{n=0}^\infty \alpha_n\xi^n.
\]
This can be integrated term by term to
\[
\bar{F}=\sum_{n=0}^\infty \beta_n\bar{\xi}^n+\alpha_n\xi^{n+1}\left({\textstyle{\frac{1}{n+1}}}
     \pm{\textstyle{\frac{2}{n+2}}}\xi\bar{\xi}+{\textstyle{\frac{1}{n+3}}}\xi^2\bar{\xi}^2\right),
\]
for complex constants $\beta_n$.
Now we impose the Lagrangian condition, that
\begin{align}
(1\pm\xi\bar{\xi})\bar{\partial}\bar{F}\mp2\xi\bar{F}=&\sum_{n=0}^\infty \left[n\beta_n\bar{\xi}^{n-1}\pm(n-2)\beta_n\xi\bar{\xi}^n
  \mp{\textstyle{\frac{2}{(n+1)(n+2)}}}\alpha_n\xi^{n+2}\right.\nonumber\\
&\left.\qquad\qquad-{\textstyle{\frac{2}{(n+2)(n+3)}}}\alpha_n\xi^{n+3}\bar{\xi}\right],\nonumber
\end{align}
is real. By a translation we can set $\beta_0=\beta_1=\beta_2=0$ and then 
${\scriptstyle{(n+1)(n+2)(n+3)}}\beta_{n+3}=\mp2\bar{\alpha}_n$ for $n\ge0$.
Letting $\alpha_n=\mp{\scriptstyle{(n+1)(n+2)(n+3)}}\lambda_n$ gives the stated result for $F$.

Finally, it is easily checked that the expressions for $r$ and $F$ satisfy (\ref{e:supfunc2}).
\end{pf}

Given a surface S in ${\Bbb{E}}^3$ (${\Bbb{E}}^3_1$) the eigen-directions of the second
fundamental form determine a pair of mutually orthogonal foliations,
called the principal foliations of S. The foliations have
singularities at points where the eigenvalues are equal,
i.e. at umbilical points (non-zero eigenvalues) or at flat points
(zero eigenvalues). By definition, the mean curvature vanishes on a
minimal (maximal) surface so the curvature is non-positive. Thus the only
singularities that can arise for the principal foliations are flat points.

\begin{Cor}
Umbilic points on minimal (maximal) surfaces are isolated and the index of the principal foliation about an 
umbilic point on a minimal (maximal) surface is less than or equal to zero.
\end{Cor}
\begin{pf}
By Proposition \ref{p:2ndfund} an umbilic point is a point where $\partial\bar{F}=0$. 
Moreover, the argument of $\bar{\partial}F$ gives the principal foliation of
the surface. Given that minimality (maximality) implies the holomorphic condition (\ref{e:mineq}), the zeros
of $\partial\bar{F}$ are isolated and have index greater than or equal to zero.
\end{pf}

The classical Weierstrass representation can be considered as the construction of a minimal surface in ${\Bbb{E}}^3$
from a holomorphic curve in ${\Bbb{L}}^{3}$ \cite{hitch}. In the following, we extend this perspective to maximal surfaces 
in ${\Bbb{E}}^3_1$:

\begin{Thm}
The surface in ${\Bbb{E}}^3$ (${\Bbb{E}}^3_1$) determined by a local holomorphic section of $TN\rightarrow N$,
$\xi\mapsto (\xi,w(\xi))$ for N=$S^2$ (H$^2$) given by
\[
z=\mp({\textstyle{\frac{1}{2}}}\xi^2\partial\partial w-\xi\partial w+w)
     +{\textstyle{\frac{1}{2}}}\bar{\partial}\bar{\partial}\bar{w}
\]
\[
t=\mp\left({\textstyle{\frac{1}{2}}}\xi\partial\partial w-{\textstyle{\frac{1}{2}}}\partial w
  +{\textstyle{\frac{1}{2}}}\bar{\xi} \bar{\partial}\bar{\partial}\bar{w}
       -{\textstyle{\frac{1}{2}}}\bar{\partial}\bar{w}\right),
\]
is minimal (maximal). Here the upper (lower) sign is the Euclidean (Lorentzian) case.
\end{Thm}
\begin{pf}
We have the push forward of the coordinate vectors:
\[
\frac{\partial}{\partial \xi}={\textstyle{\frac{1}{2}}}\partial\partial\partial w
   \left(\mp\xi^2\frac{\partial}{\partial z}+\frac{\partial}{\partial \bar{z}}
    \mp\xi\frac{\partial}{\partial t}\right).
\]
The unit vector which corresponds to the point $\xi\in N$ is
\[
e_0=\frac{2\xi}{1\pm\xi\bar{\xi}}\frac{\partial}{\partial z}
   +\frac{2\bar{\xi}}{1\pm\xi\bar{\xi}}\frac{\partial}{\partial \bar{z}}
   +\frac{1\mp\xi\bar{\xi}}{1\pm\xi\bar{\xi}}\frac{\partial}{\partial t}.
\]
The inner product of the preceding 2 vectors is zero, and so $\xi$ is the unit normal direction to the surface.

The relationship with the fibre coordinate follows from equation (\ref{e:coord2})
\[
\eta={\textstyle{\frac{1}{2}}}\left(z-2t\xi\mp\bar{z}\xi^2\right),
\]
and this turns out to be
\[
\eta={\textstyle{\frac{1}{4}}}(1\pm\xi\bar{\xi})^3\bar{\partial}\bar{\partial}  
   \left(\frac{\bar{w}}{1\pm\xi\bar{\xi}}\right)\mp{\textstyle{\frac{1}{2}}}w.
\]
Finally, differentiating again we find that
\[
\frac{\bar{\partial} \eta}{(1\pm\xi\bar{\xi})^2} 
    ={\textstyle{\frac{1}{4}}}\bar{\partial}\bar{\partial}\bar{\partial}\bar{w} .
\] 
Thus, by the identity (\ref{e:cm3}), $r+\rho_0=C$, and the line congruence is orthogonal
to a minimal (maximal) surface. That this surface is in fact given parametrically as claimed can be 
seen by inserting these in equation (\ref{e:coord2}) and computing that
\[
r+(1\pm\xi\bar{\xi})^2\partial\left(\frac{\eta}{(1\pm\xi\bar{\xi})^2}\right)=0.
\]
\end{pf}

\noindent{\bf Acknowledgement}:

The authors would like to thank Marcos Salvai for helpful discussions.

This work was supported by the Research in Pairs Programme of the Mathematisches Forschungsinstitut Oberwolfach, Germany.


\begin{thebibliography}{10}

\bibitem{gak1}
B. Guilfoyle and W. Klingenberg, {\it On the space of oriented affine
  lines in ${\Bbb{E}}^3$}, Archiv der Math. {\bf 82} (2004), 81--84.

\bibitem{gak2}
B. Guilfoyle and W. Klingenberg, {\it Generalised surfaces in
  ${\Bbb{R}}^3$}, Math. Proc. of the R.I.A. {\bf 104A} (2004) 199--209.

\bibitem{gak4}
B. Guilfoyle and W. Klingenberg, {\it An indefinite K\"ahler metric on the space of oriented lines},  J. London Math.
Soc. {\bf 72}, (2005) 497--509.

\bibitem{hitch}
N.J. Hitchin, {\it Monopoles and geodesics}, Comm. Math. Phys. {\bf 83} (1982) 579--602. 

\bibitem{kob}
O. Kobayashi, {\it Maximal surfaces in 3-dimensional Minkowski space ${\Bbb{L}}^3$}, Tokyo J. Math. {\bf 6} (1983) 
297--309.


\bibitem{KN}
S. Kobayashi and K. Nomizu, Foundations of Differential Geometry, Wiley, New York (1996).

\bibitem{salvai}
M. Salvai, {\it On the geometry of the space of oriented lines in Euclidean space}, Manuscripta Math. {\bf 118} (2005)
181--189.

\end{thebibliography}
\end{document}